\documentclass{amsart}
 
\usepackage{amssymb, hyperref}
 
\usepackage{ amsmath,amssymb,amsbsy,amsfonts, latexsym,amsopn,amstext, amsxtra,euscript,amscd}

\usepackage{amssymb}
\usepackage{latexsym} 
\usepackage{url}

\usepackage{hyperref}

\usepackage{amsrefs}

\newtheorem{thm}{Theorem}
\newtheorem{lem}{Lemma}

\def\Z{\mathbb Z}

\def\Q{\mathbb Q}

\def\P{\mathbb P}

\def\X{\mathcal X}

\def\O{\mathcal O}
\def\p{\mathfrak p}

\def\v{\mathfrak v}
\newcommand\m{\mathfrak m}

\def\s{\mathfrak s}

\def\X{\mathcal X}

\def\d{{\delta }}
\def\a{{\alpha }}

\def\D{\Delta}

\def\D{\Delta}

\def\l{\lambda}
\def\a{\alpha}

\begin{document}


\title{Equations of curves with minimal discriminant}

\author{Rachel Shaska}
%

\address{Department of Electrical Engineering, Oakland University, Rochester, MI, 48386.}

\email{rishaska@oakland.edu}
%






\def\dis{\mathfrak D}
\def\s{\sigma}
\def\<{\left<}
\def\>{\right>}

\begin{abstract}
In this paper we give an algorithm of how to determine a Weierstrass equation with minimal discriminant for superelliptic curves generalizing work of Tate \cite{ta-75} for elliptic curves and Liu \cite{liu} for genus 2 curves.
\end{abstract}

\maketitle

\section{Introduction}
%
Let $K$ be a field with a discrete valuation $\v$  and ring of integers $\O_K$   and $C$ an irreducible, smooth, algebraic    curve of genus $g\geq 1$ defined over $K$ and function field $K (C)$. The discriminant $\dis_{C/K}$ is an important invariant of the function field of the curve and therefore of the curve. Since the discriminant is a polynomial given in terms of the coefficients of the curve, then it is an ideal in the ring of integers $\O_K$ of $K$.  The valuation of this ideal is a positive integer.  
A classical question is to find an equation of the curve such that this valuation is minimal, in other words the discriminant is minimal.  

When $g=1$, so that $C$ is an elliptic curve, there is an extensive theory  of the minimal discriminant ideal $\dis_{C/K}$.  Tate \cite{ta-75}  devised an algorithm how to determine the Weierstrass equation of an elliptic curve with minimal discriminant as part of his larger project of determining Neron models for elliptic curves.  The main focus of this paper is to extend their work to  superelliptic curves, full details and proofs are intended in \cite{rachel}.

The paper is organized as follows. In section 2 we give the basic definitions for genus $g\geq 2$ superelliptic curves isomorphism classes of which correspond to projectively equivalent classes of degree $d$ binary forms. For a  binary form $f(X, Z)$ and a matrix 
$M =  \begin{bmatrix} a, b \\ c, d \end{bmatrix}$, such that  $M \in GL_2 (k)$, we have that $f^M:= f(aX + bZ, cX + dZ)$ has discriminant 
$\D (f^M) = ( \det M )^{d(d-1)} \cdot \D (f)$. 
This property of the discriminant is crucial in our algorithm which is explained in Section~5.

In section 3 we define the discriminant of a genus $g \geq 2$ superelliptic curve $\X_g$ defined over an algebraically number field $K$. We follow the classical theory and define the discriminant for local fields and then generalize it to global fields. 

In Section 4, we summarize briefly Tate's algorithm and a modified version of it by Laska \cite{la-82}. Since the case of the elliptic curves is the simplest case this hopefully gives the reader an idea of how things work out in higher genus.  
In Section 5, we generalize the algorithm to all superelliptic curves.  This algorithm  computes a Weierstrass equation with minimal discriminant for all superelliptic curves.  Details and  proofs are intended to be described in \cite{rachel}.

\section{Preliminaries}

Let $\X_g$ be a superelliptic curve of genus $g\geq 2$   with affine equation 
\begin{equation}\label{w-eq-super}
y^n = f(x, 1) = a_d x^d + \cdots a_1 x + a_0
\end{equation}
defined over and algebraic number field $K$.   Obviously the set of roots of $f(x)$ does not determine uniquely the isomorphism class of $\X_g$ since every coordinate change in $x$ would change the set of these roots.  Such isomorphism classes are classified by the invariants of binary forms.   

%
For any algebraically closed field $k$ let  $k[X,Z]$  be the  polynomial ring in  two variables and  let $V_d$ denote  the $(d+1)$-dimensional  subspace  of  $k[X,Z]$  consisting of homogeneous polynomials
\begin{equation}  \label{eq1}
f(X,Z) = a_0X^d + a_1X^{d-1}Z + ... + a_dZ^d
\end{equation}
of  degree $d$. Elements  in $V_d$  are called  \textit{binary  forms} of degree $d$.    $GL_2(k)$ act as a group of automorphisms on $ k[X, Z] $   as follows:
\begin{equation}
 M =
\begin{pmatrix} a &b \\  c & d
\end{pmatrix}
\in GL_2(k), \textit{   then       }
\quad  M  \begin{pmatrix} X\\ Z \end{pmatrix} =
\begin{pmatrix} aX+bZ\\ cX+dZ \end{pmatrix}
\end{equation}
Denote by $f^M$ the binary form $f^M (X, Z) := f(aX+bZ, cX+dZ)$. 
It is well  known that $SL_2(k)$ leaves a bilinear  form (unique up to scalar multiples) on $V_d$ invariant. 

Consider $a_0$, $a_1$,  ... , $a_d$ as parameters  (coordinate  functions on $V_d$). Then the coordinate  ring of $V_d$ can be identified with $ k[a_0  , ... , a_d] $. For $I \in k[a_0, ... , a_d]$ and $M \in GL_2(k)$, define $I^M \in k[a_0, \dots , a_d]$ as follows
\begin{equation} \label{eq_I}
{I^M}(f):= I( f^M)
\end{equation}
for all $f \in V_d$. Then  $I^{MN} = (I^{M})^{N}$ and Eq.~(\ref{eq_I}) defines an action of $GL_2(k)$ on $k[a_0, \dots , a_d]$.
A homogeneous polynomial $I\in k[a_0, \dots , a_d, X, Z]$ is called a \textbf{covariant}  of index $s$ if  $I^M(f)=\delta^s I(f)$, 
where $\delta =\det(M)$.  The homogeneous degree in $a_0, \dots , a_d$ is called the \textbf{ degree} of $I$,  and the homogeneous degree in $X, Z$ is called the \textbf{  order} of $I$.  A covariant of order zero is called \textbf{invariant}.  An invariant is a $SL_2(k)$-invariant on $V_d$.

Let $f(X, Z)$ and $g(X, Z)$   be binary forms of  degree $n$ and $m$ respectively with coefficients in $k$. We denote the \textbf{r-transvection} of two binary forms $f$ and $g$ by $(f,g)^r$.
It is a homogeneous polynomial in $k[X, Z]$ and therefore a covariant of order $m+n-2r$ and degree 2.  

%

A very important  invariant is the discriminant of the binary form.  In the classical way, the discriminant is defined as 
$ \D =  \prod_{i \neq j } (\a_i - \a_j)^2$, 
where $\a_1, \dots \a_d$ are the roots of $f(x, 1)$.   It is a well-known result that it can be expressed in terms of the transvectians. For example, for binary sectics we have $\D = J_{10}$  and for binary octavics    $\D (f) = J_{14}$; see \cite{rachel} for details. 

\begin{lem}
i) The discriminant of a degree $d$ binary form $f(X, Z)\in k[X, Z]$ is and $SL_2 (k)$-invariant of degree $2d-2 $.  

ii) For any $M \in GL_2 (k)$ and any degree $d$ binary form $f$ we have that 
\[ \D (f^M) = \left(  \det M \right)^{d (d-1) } \, \D (f) \]
\end{lem}
%

\section{Discriminant of a curve}
The concept of a minimal discriminant for elliptic curves was defined by Tate and others in the 1970-s; see \cite{ta-75}.    Such definitions and results we generalized by Lockhart in \cite{lockhart} for hyperelliptic curves.  In this section we briefly generalize the concept of  the minimal discriminant to all superelliptic curves. 

%
Let $K$ be a local field, complete with respect to a valuation $\v$.  Let $\O_K$ be the ring of integers of $K$, in other words 
$\O_K = \{ x\in K \, | \, \v (x) \geq 0\}$. 
We denote by $\O_K^\ast$ the group of units of $\O_K$ and by $\m$ the maximal ideal of $\O_K$.  Let $\pi$ be a generator  for $\m$ and $k=\O_K / \m$ the residue field. We assume that $k$ is perfect and denote its algebraic closure by $\bar k$. 

Let $\X_g$ be a superelliptic curve of genus $g \geq 2$ defined over $K$ and $P$ a $K$-rational point on $\X_g$.     By a suitable change of coordinates we can assume that all coefficients of $\X_g$ are in $\O_K$.      Then, the discriminant $\D \in \O_K$.  In this case we say that the equation of $\X_g$ is \textbf{integral}.

An equation for $\X_g$ is said to be a \textbf{minimal equation}  if it is integral and $\v (\D)$ is minimal among all integral equations of $\X_g$. The ideal $I=\m^{\v (\D)}$ is called the \textbf{minimal discriminant} of $\X_g$. 
 
\newcommand\fa{\mathfrak a}

%
Let us assume now that $K$ is an algebraic number field with field of integers $\O_K$.  Let $M_K$ be the set of all inequivalent absolute values on $K$  and $M_K^0$ the set of all non-archimedean absolute values in $M_K$. 
We denote by $K_\v$ te completion of $K$ for each $\v \in M_K^0$ and by $\O_\v$ the valuation ring in $K_\v$. Let $\p_v$ be the prime ideal in $\O_K$ and $\m_v$ the corresponding maximal ideal in $K_\v$. Let $(\X, P)$ be a superelliptic curve of genus $g\geq 2$ over $K$.  

If $\v \in M_K^0$ we say that $\X$ is \textbf{integral at $\v$} if $\X$ is integral when viewed as a curve over $K_\v$.  We say that $\X$ is \textbf{minimal at $\v$} when it is minimal over $K_\v$. 

An equation of $\X$ over $K$ is called \textbf{integral} (resp. \textbf{minimal}) over $K$ if it is integral (resp. minimal) over $K_\v$, for each $\v \in M_K^0$. 
 
Next we will define the minimal discriminant over $K$ to be the product of all the local minimal discriminants. For each $\v \in M_K^0$ we denote by $\D_\v$ the minimal discriminant for $(\X, P)$ over $K_\v$.  The \textbf{minimal discriminant} of $(\X, P)$ over $K$ is the ideal 
\[ \D_{\X / K} = \prod_{\v \in M_K^0} \m_\v^{\v (\D_\v) } \]
We denote by $\fa_\X$ the ideal   $ \fa_\X = \prod_{\v \in M_K^0} \p_\v^{\v (\D_\v) }$. 
In \cite{rachel} we prove that 
\begin{thm}
Let $(\X_g, P)$ be a superelliptic curve  over $\Q$. Then its global minimal  discriminant $\D\in \Z$  is unique (up to multiplication by a unit).  There exists a minimal Weierstrass equation corresponding to this $\D$.
\end{thm}
Next we briefly describe how this minimal Weierstrass equation is determined for superelliptic curves.  Full details and further analysis of discriminants of superelliptic curves is intended in \cite{rachel}. 
 

\section{Elliptic curves and Tate's algorithm}
Let $E$  be an elliptic curve defined over a number field $K$ with equation
\begin{equation}\label{w-eq-ell}     y^2 + a_1 xy + a_3  y = x^3 + a_2 x^2 + a_4  x +a_6. \end{equation}
For simplicity we assume that $E$ is defined over $\Q$, the algorithm works exactly the same for   any algebraic number field $K$.

We would like to find an equation
\begin{equation}\label{w-eq-2}     y^2 + a_1^\prime xy + a_3^\prime y = x^3 + a_2^\prime x^2 + a_4^\prime x +a_6^\prime. \end{equation}
such that the discriminant $\D^\prime$ of the curve in Eq.~\eqref{w-eq-2} is minimal.  Since we want the new equation to have integer coefficients then the only transformations we can have are 
\[ x = u^2 x^\prime + r, \qquad y = u^3 y^\prime + u^2 s x^\prime + t\]
for $u, r, s, t \in \Z$ and $u \neq 0$.  The coefficients of the two equations are related as follows:
\[
\begin{split} 
& ua_1^\prime  = a_1 + 2s,                          \\
& u^3 a_3^\prime  = a_3 + r a_1 + 2 t,               \\          
& u^2 a_2^\prime  = a_2 - sa_1 + 3 r - s^2,             \\
\end{split}
\qquad
\begin{split}
                   & u^4 a_4^\prime  = a_4 - s a_3 + 2 r a_2 - (t+rs)a_1 + 3 r^2 - 2st    \\
      & u^6 a_6^\prime  = a_6 + r a_4 + r^2 a_2 + r^3- ta_3 - rta_1 - t^2   \\          
    & u^{12} \D^\prime  =  \D                                              \\
\end{split}
\]
The version of the algorithm below is due to M. Laska; see \cite{la-82}. \\

\noindent \textsc{Step 1:} Compute the following
\[
\begin{split}
 c_4 & = (a_1^2+4a_2)^2 - 24(a_1a_3 +2a_4), \\
  c_6  & = - (a_1^2+4a_2)^3 + 36(a_1^2 + 4 a_2)(a_1a_3 + 2 a_4) - 216 (a_3^2 +4a_6) 
\end{split}
\]
\noindent  \textsc{Step 2:} Determine the set $S$ of integers $u \in \Z $ such that there exist $x_u$, $y_u \in \Z$ such that 
$ u^4 = x_u c_4$ and $ u^6 y_u = c_6$.  Notice that   $S$ is a finite set. \\

\noindent  \textsc{Step 3:} Choose the largest $u \in S$, say $u_0$ and factor it as $u_0 = 2^{e_2} \, 3^{e_3} \, v$, where $v$ is relatively prime to 6.  \\

\noindent  \textsc{Step 4:}  Choose 
\[
a_1^\prime, a_3^\prime \in \left\{ \sum_{i=1}^n \alpha_i w_i \, | \, \alpha_i = 0 \, \textit{ or }  1 \, \right\} \, \textit{ and } \,
a_2^\prime \in \left\{ \sum_{i=1}^n \alpha_i w_i \, | \, \alpha_i = -1, 0 \, \textit{ or }  1 \, \right\}
\]
subject to the following conditions:
\[ (a_1^\prime)^4 \equiv x_u \mod 8, \quad (a_2^\prime)^3 \equiv - (a_1^\prime)^6 - y_u \mod 3.  \]

\noindent  \textsc{Step 5:}  Solve the following equations for $a_4^\prime$ and $a_6^\prime$
\[
\begin{split}
 x_u & = ({a_1^\prime}^2+4{a_2^\prime})^2 - 24({a_1^\prime}{a_3^\prime} +2{a_4^\prime}), \\
 y_u  & = - ({a_1^\prime}^2+4{a_2^\prime})^3 + 36({a_1^\prime}^2 + 4 {a_2^\prime})({a_1^\prime}{a_3^\prime} + 2 {a_4^\prime}) - 216 ({a_3^\prime}^2 +4{a_6^\prime}) 
\end{split}
\] 

\noindent  \textsc{Step 6:} Solve the equations for $s, r, t$ successively 
\[ u a_1^\prime = a_1 + 2s, \quad u^2 a_2^\prime = a_2 - s a_1 + 3 r - s^2, \quad u^3 a_3^\prime = a_3 + r a_1 + 2t \]
For these values of $a_1^\prime, \dots , a_6^\prime$  the Eq.~\eqref{w-eq-2}  is the desired result.

For a complete version of the algorithm see \cite{la-82}.  
 
\section{Superelliptic curves with minimal discriminant}
Let $\X_g$ be a genus $g\geq 2$ superelliptic curve with equation as in \eqref{w-eq-super}. 
The discriminant of $\X_g$ is the discriminant of the binary form $f(x, z)$, hence an invariant of homogenous degree $\d= 2d -2$ and $\D_f \in \O_K$. 

Let $M \in GL_2 (K)$ such that $\det M = \l$. Then from remarks in section 2 we have that $ \D (f^M) =   \l^{d(d-1)} \,  \D(f)$. 
We perform the coordinate change $x \to \frac 1 {u^n} x$ on  $f(x)$.  Then the new discriminant is  $\D^\prime = \frac 1 {u^{n \cdot d(d-1)} } \, \cdot  \D $.

\begin{lem}  A superelliptic curve $\X_g$ with integral equation  
\[ y^n = a_d x^d + \cdots a_1 x + a_0 \]
 is in minimal form if $\v (\D) <  n d(d-1)$.
\end{lem}

Hence, if we choose $u\in \Z$ such that $u^{ n d(d-1) }$  divides $\D$, then $\D^\prime$ becomes smaller.  Indeed, we would like to choose the largest such $u$.  In the process we have to make sure that for the $u$'s that we pick we do get an equation of a superelliptic curve isomorphic to $\X_g$. 

Hence, we factor $\D$ as a product of primes, say $\D = p_1^{\alpha_1} \cdots p_r^{\alpha_r}$,   and take $u$ to be the product of those powers of primes with exponents $\alpha_i \geq n d (d-1)$.  For primes $p=2, 3$ we have to be more careful since in our exposition above we have assumed that the characteristic of the field is $\neq 2, 3$. 

In \cite{rachel}  we give a description of all the steps of the algorithm which is also implemented for genus 3 hyperelliptic curves and for triagonal curves   $y^3=f(x)$ up to $\deg f \leq 8$.  

The main result of \cite{rachel} is the following:

\begin{thm}
Let $\X$ be a genus $g \geq 2$ superelliptic curve defined  over an algebraic number field $K$ and $P$ a $K$-rational point on $\X$. For the pair $(\X, P)$  the
 global minimal  discriminant $\D_{\X, P}   \in \O_K$  is unique (up to multiplication by a unit).  Moreover, there exists a minimal Weierstrass equation corresponding to this discriminant $\D_{\X, P}$.
\end{thm}

The theorem also provides the blueprint for the algorithm which is much more involved then the case of elliptic curves described in \cite{ta-75} and \cite{la-82} and the case of hyperelliptic curves described in \cite{liu}.

\bibliographystyle{amsplain}

\end{document}